# Comfort-constrained Distributed Heat Pump Management

Simon Parkinson, Dan Wang, Curran Crawford, and Ned Djilali

*Abstract*— This paper introduces the design of a demand response network control strategy aimed at thermostatically controlled electric heating and cooling systems in buildings. The method relies on the use of programmable communicating thermostats, which are able to provide important component-level state variables to a system-level central controller. This information can be used to build power density distribution functions for the aggregate heat pump load. These functions lay out the fundamental basis for the methodology by allowing for consideration of customer-level constraints within the system-level decision making process. The proposed strategy is then implemented in a computational model to simulate a distribution of buildings, where the aggregate heat pump load is managed to provide the regulation services needed to successfully integrate wind power generators. Increased exploitation of wind resources will place similarly themed ancillary services in high-demand, traditionally provided by dispatchable energy resources that are ill-suited for the frequent power gradients that accompany wind power generation.

*Keywords*- Demand Response, Heat Pumps, Programmable Communicating Thermostats, Wind Energy Integration

## I. Introduction

With the increasing prevalence of information technology and wireless networking capabilities, power systems can interact with the loads they service, inducing responsive devices into operational trajectories that are beneficial from a system-level perspective. The main drawback lies in the need to balance these objectives with those expected by the customer at end-use. This paper addresses this issue, and provides a control strategy capable of accurately controlling distributed heat pump systems, constrained by the comfort-level commensurate with customer satisfaction. The primary power system objective is to manage these units to provide grid-side services that routinely result in ineffective utilization of conventional energy resources.

Heat pumps represent an extremely efficient method of providing both heating and cooling functions within a single unit. As these devices rely upon the electrical grid for energy, emissions related to conventional methods of heating (natural gas- or oil-based furnaces) can be avoided if these grids employ renewable energy resources [1]. Many of these resources, such as wind and wave energy, display considerable short-term variability, making integration into legacy power systems difficult, and would in turn benefit greatly from demand that displays flexibility [2]. Furthermore, the increased capacity requirements brought on by shifting heating demand to the electrical load can be met efficiently with co-located renewable generation, as transmission losses would be avoided. Thus, this work focuses on applying the proposed distributed heat pump control strategy in a local energy system incorporating wind power generators.

## II. Formulation of the Management System

### A. Component-level Control

To obtain fully responsive non-disruptive control, the method first proposed by Callaway in [3] is employed. Under this strategy, partial synchronization of thermostatically controlled loads provides the dynamic response in power demand through small changes $u$ to the temperature set-point $\theta_s$. Each heat pump consumes power based upon the indoor air temperature $\theta_a$ measured at the conclusion of interval $k$. The machine-state $n$ (equal to 1 for active, and 0 for inactive) of population element $i$ can be captured for a heating process as follows:

$$n_i(k+1) = \begin{cases} 1 & \theta_{a,i}(k) \leq \theta_- = \theta_s - \frac{\delta}{2} + u(k) \\ 0 & \theta_{a,i}(k) \geq \theta_+ = \theta_s + \frac{\delta}{2} + u(k) \\ n_i(k) & otherwise \end{cases} \quad (1)$$

where $\theta_-$ and $\theta_+$ describe the lower and upper boundaries of the temperature deadband that spans the distance $\delta$. This deadband exists to prevent any rapid cycling of the device that would accompany a definite set-point objective.

Programmable communicating thermostats (PCTs) represent an innovative smart-grid technology that allows for networked communication with outside entities. Routinely these devices drive load control strategies through consideration of dynamic pricing schedules [4], however in these scenarios problems related to response stability have been uncovered [5]. For this reason, in this work these devices are used to directly communicate pertinent state variables associated with each building-heat pump system, which in this case is defined as the *power-state vector* of

Manuscript received March 19, 2011, revised June 29, 2011. This work was supported in part by the Pacific Institute for Climate Solutions, NSERC-H2CAN, and NSERC-WESNet.

S. Parkinson, D. Wang, C. Crawford, and N. Djilali are with the Dept. of Mech. Eng., and the Institute for Integrated Energy Systems, University of Victoria, PO Box 3055 STN CSC, Victoria BC V8W 3P6, Canada. (e-mail:scp@uvic.ca).

element *i*, given by:

$$\mathbf{x_i}(k) = [\, n_i(k) \;\; \theta_{a,i}(k) \;\; P_{h,i} \,]^T \qquad (2)$$

where $P_{h,i}$ is the heat pump's rated power consumption when active. In our system configuration, once a PCT has communicated its power-state vector, it waits before reacting to the control algorithm given by (1), thereby allowing for analyses based upon current component-level inputs. The resultant delay-time must be considerably less than the thermostat update time (approximately one minute) to ensure compatibility, and will be directly related to the number of participating customers, as well as the efficiency of system-level decision making.

## B. System-level Control

The power-state vectors communicated to the central controller can be used to generate *power density distribution functions* (PDDF) for both the active ($\phi_1$) and inactive ($\phi_0$) machine-states. These functions describe the amount of power at a given indoor air temperature relative to the total installed power that exists in the population $P_{H,cap}$. This capacity rating corresponds to the sum of the individual device ratings (third row of each power-state vector). The total aggregate heat pump load is then $P_{H,cap}$ multiplied by the power density in the active-state, which we will define as the *capacity-factor* $\Phi$ of the heat pump load.

$$P_H(t) = \sum_{i=1}^{N} P_{h,i} \int_{-\infty}^{\infty} \phi_1(t,\theta)\,d\theta = P_{H,cap}\Phi(t) \qquad (3)$$

As each element in the population waits for a central response upon sending its power-state vector, all are basically synchronized to the same clock. This means that the discrete thermostatic control event given by (1) introduces a discontinuity between sampling intervals, wherein the air temperature measurement is utilized to determine the individual machine-states thereafter, thus re-distributing the power density distribution in the aggregate system accordingly. Any distribution in a given state that has traversed past the corresponding state-transition boundary ($\theta_+$ for *n* = 1, and $\theta_-$ for *n* = 0), will be transferred to the opposite distribution. System-level objectives will need to focus on controlling this process, as the amount of distribution in the active-state dictates the aggregate heat pump load. Boundary conditions on either side of the discontinuity event require that:

$$\int_{-\infty}^{\infty} \phi_1(k+1,\theta)\,d\theta = \int_{-\infty}^{\theta_-} \phi_0(k,\theta)\,d\theta + \int_{-\infty}^{\theta_+} \phi_1(k,\theta)\,d\theta \qquad (4)$$

where the righthand-side PDDFs are strictly defined by the power-state vectors communicated to the central controller. Examining (3), it is clear that although (4) represents the state-equation directly relating control-input (set-point temperature) to the capacity-factor, this expression will be difficult to utilize in the current integral form.

Constraining the thermostats to a certain integer resolution *R* over a temperature range $\theta_R$ will break the possible measurements that can occur into discrete intervals $\Delta\theta$. For the considered system $\theta_R$ is taken to be twice the deadband width.

$$\Delta\theta = \frac{\theta_R}{R} = \frac{2\delta}{R} \;;\; R \in \mathbb{N} \qquad (5)$$

The temperature measurements that can occur are thus:

$$\theta(m) = \theta_s - \delta\left(1 - \frac{2m}{R}\right) \;;\; m = 0, 1, 2, ..., R \qquad (6)$$

Upon discretizing (4), a *capacity-factor function* (CFF) with respect to set-point index $m_s$ can be obtained:

$$\Phi(k+1, m_s) = \sum_{m=0}^{\epsilon_-} \phi_0(k,m)\Delta\theta + \sum_{m=0}^{\epsilon_+} \phi_1(k,m)\Delta\theta$$
$$\epsilon_+ = m_s + \frac{\delta}{2\Delta\theta} \;;\; \epsilon_- = m_s - \frac{\delta}{2\Delta\theta} \qquad (7)$$

which for component-level thermostats operating under the equivalent constrained measurement scenario is exact.

Accurately scheduling the heat pump load in real-time will initially require determination of the control signal (set-point change) corresponding to the objective. To prevent customer-side disruption, set-point modulations are constrained to remain within the quarter-deadband width, as it is unlikely individuals will notice changes of this magnitude at end-use.

$$|\,u(k)\,| \leq \frac{\delta}{4} \qquad (8)$$

Constraining set-point changes to this magnitude provides further benefits, as it will only involve loads traversing the final quarter-trajectory of their current operating state, thereby preventing rapid-cycling. By utilizing (6) to compute (7) at these limits, the current CFF obtained can be used to determine the constraints given by (8) at the system-level:

$$\Phi(k+1, m_{s,min}) \leq \Phi(k+1, m_s) \leq \Phi(k+1, m_{s,max}) \qquad (9)$$

where $m_{s,min}$ is the minimum feasible set-point index, with $m_{s,max}$ denoting the maximum.

The feasible region given by (9) limits the possible power-gradients which the aggregate heat pump load can withstand without compromising end-service. This interval directly relates customer-level constraints (thermal comfort the set-point temperature provides) to system-level planning (aggregate operating capacity of the heat pump load), and can now be utilized in defining the desired capacity factor $\Phi^*$ that should be seen from the aggregate system (the tracking signal). As the PDDFs are defined to be:

$$\phi_0(k,m)\,,\; \phi_1(k,m) \geq 0 \;\; \forall \;\; k, m \qquad (10)$$

then from (7), the CFF is necessarily increasing with respect to $m_s$, meaning the following deviation minimization

problem:

$$\text{Min} \quad F = [\,\Phi^*(k+1) - \Phi(k+1, m_s)\,]^2 \quad (11)$$

is strictly convex, guaranteed only one critical point $m_s^*$ over the interval bounded by (9). Once $\Phi^*$ is defined, solution to the simple minimization problem can be obtained, and used to calculate the optimal set-point change.

$$u(k) = \theta\,(m_s^*) - \theta_s = \delta\left(\frac{2m_s^*}{R} - 1\right) \quad (12)$$

This value is integrated into the component-level thermostats through communication with the individual PCTs, which will then respond according to (1).

## III. Applications

### A. Signal Tracking

An example of the proposed method is given in Fig. 1, where a heterogeneous group of 1000 building-heat pump transient models, of similar structure to [6], are coupled to thermostats that update each interval based on (1). The tracking signal $\Phi^*$ is applied for $k \geq 100$, and has been arbitrarily defined as the steady-state level perturbed by a stochastic disturbance, constrained to remain within the feasible region defined by (9). A static comfort condition of 20°C and deadband width of 1°C has been used for each building, with a thermostat resolution of $R = 1000$ pursued to accelerate computations. At the conclusion of each sampling interval, the PDDFs are determined from the power-state vectors obtained from evolving the building models. The feasible region is then enumerated, and a target capacity-factor defined. The optimal set-point index is then computed using (11), with the control signal sent to the component-level models then defined by (12). The result is the physical aggregate load trajectory corresponding to the target-level over each interval applied.

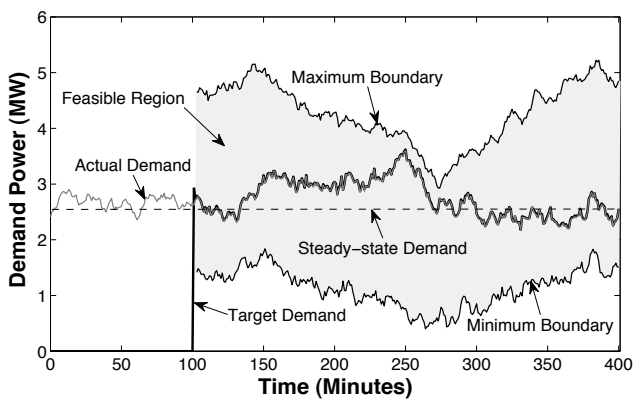

Fig. 1. Control of 1000 heat pump systems under the proposed methodology, where the actual demand is determined by simulating transient models coupled to the thermostats given by (1). A static outdoor temperature of 4°C has been used.

Between the region ($150 \leq k \leq 250$), the desired trajectory demands power-levels greater than the steady-state rating over a range much greater than the average time it takes individual systems to traverse their entire operating cycle (approximately 45 min). As the desired aggregate heat pump demand is greater than the rate at which the system is dissipating energy to the surroundings (the steady-state losses), the amount of energy stored within the buildings must increase. Thus, all buildings effectively *charge* to their limit, pushing the indoor air temperatures nearer to the upper state-transition boundary. As the thermostats at the customer-level are constrained by (8) and (1), $\Phi^*$ must decrease to remain in control. By tracking the PDDFs, stability is guaranteed through consideration of (9) in defining the target, which in Fig. 1 is seen to shift the feasible region to reflect the customer's comfort constraints.

### B. Regulating Wind Power Variability

Next, the proposed management strategy is utilized to regulate power fluctuations from two 2.5 MW wind turbines integrated into a community consisting of 2000 small- to medium-scale buildings. Three types of loads are defined: a nominal load $P_N$, the heat pump load $P_H$, and the wind power input $P_W$. The nominal load describes all electrical energy demand that is not caused by the heat pumps, and is obtained by fitting curves to normalized data for typical residential houses found in [7]. Nominal load fluctuations are taken to be normally distributed, with a standard deviation equal to 10% of the off-peak levels. The wind power data are generated with a typical turbine power-curve, using wind speed data taken from [8], which has been simultaneously measured with the dynamic outdoor temperature data (average of 8°C) input to the building models that generate the heat pump load. The idealized total load is then:

$$P_L(k) = P_N(k) + P_H(k) - P_W(k) \quad (13)$$

To regulate the load we assume that the $k + 1$ nominal load and wind power output are known, and then set the target capacity-factor so as to minimize deviations from the average total load over the last two sampling intervals.

$$\Phi^*(k+1) = \frac{P_W(k+1) - P_N(k+1) + \frac{1}{2}[P_L(k) + P_L(k-1)]}{P_{H,cap}} \quad (14)$$

This perfect knowledge assumption is fair, as the control takes place at the conclusion of each thermostat update interval, meaning near real-time values for both variables could be obtained over a network similar to that controlling the heat pumps.

Results are given in Fig. 2, where two scenarios are depicted: with (controlled) and without (uncontrolled) the proposed distributed heat pump management strategy. Large load fluctuations observed in **(C)** for the uncontrolled case are ill-suited for grid-integration, and will require the use of ancillary resources for smoothing before conventional generation can meet the demand. Applying the distributed heat pump management strategy given by (14), the required leveling is achieved, with the set-point modulation computed in real-time and sent to the thermostats given in **(D)**, resulting in the controlled heat pump demand seen in **(B)**. As can be seen, these small changes to the temperature set-point achieve the control objective, with the average temperature in each building remaining near the customer's desired comfort conditions at all times. The power gradient observed in the total load can be obtained by finite-differencing the profiles. The result is given in Fig. 3 as a probability density plot over the simulated time horizon. The control strategy is seen to filter larger fluctuations, meaning building thermal mass distributed throughout the system is acting as an effective variability buffer.

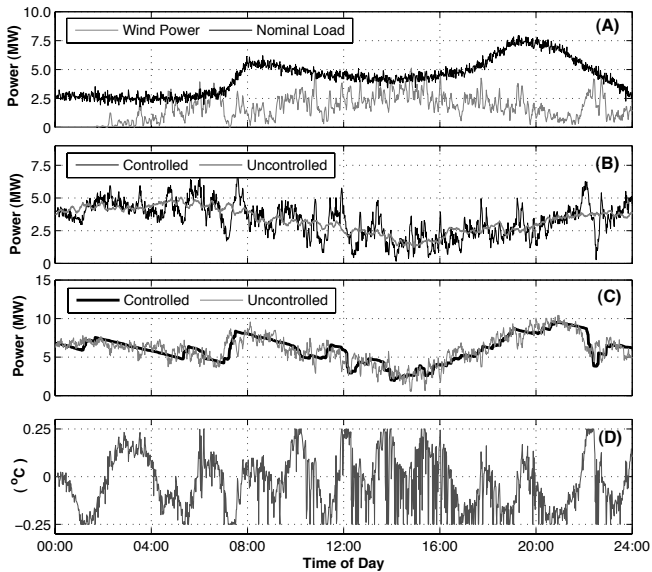

Fig. 2. Different load profiles implemented in the two scenarios considered are given in **(A)** and **(B)**, with the total load seen in **(C)** computed using (13). The control signal determined in real-time and implemented in each individual thermostat is given in **(D)**, and results in the controlled heat pump demand shown in **(B)**.

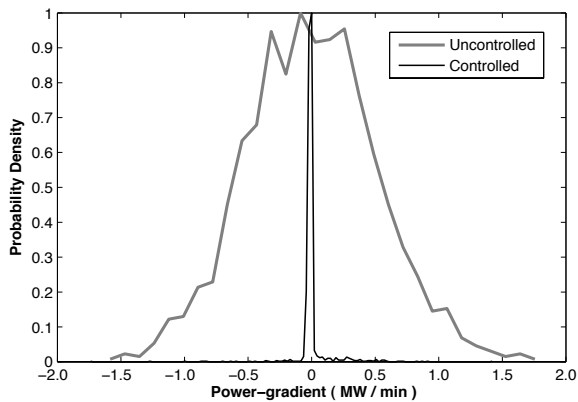

Fig. 3. Peak-normalized probability density plot of the total load power gradient for both scenarios observed in Fig. 2 **(C)**

## IV. Conclusion

This paper has introduced a new method of accurately controlling groups of electric heating and cooling systems in buildings, whilst remaining true to the comfort-levels commensurate with customer satisfaction at end-use. The main outcome of the work is the ability to consider the customer-level constraints within the system-level decision making process. The observed response displays attributes akin to frequency regulation and energy imbalance services, and as such the method is well suited for participation in similarly themed energy markets. Traditionally, these services result in inefficient operation of conventional dispatchable resources, meaning there are considerable economic and environmental benefits to be realized if these alternative methods prove successful. As the devices pursued in this work are the property of community members, these individuals look to benefit directly, through either new demand-side infrastructure or reduced costs of energy-related services. Indeed, the approach pursued in this work could enable community-owned wind energy projects, as the strategy integrates with the community-based thermal conditioning units.

Future work will focus on extending the methodology to dynamic comfort conditions. It is thought that by consideration of a third machine-state (transitionary-state), and tracking the distance from the set-point temperature instead of the indoor air temperature itself, could allow for a similar structure to be implemented. Furthermore, the spatially distributed nature of the loads has been neglected, and must be considered within distribution-level power-flow analysis. Finally, coupling of the proposed model with market-based simulations will be very beneficial to understanding the economic feasibility of actually pursuing such system configurations in grid connected applications.